# On groups whose subnormal abelian subgroups are normal

L.A. Kurdachenko, J. Otal, I.Ya. Subbotin

**Abstract:** In the current paper we study the groups, whose subnormal abelian subgroups are normal. We obtained a quite detailed description of such hyperabelian groups with a periodic Baer radical. The description of hyperabelian Lie algebras, whose abelian subideals are ideals, is also obtained.

Key words: subnormal abelian subgroups, hyperabelian groups, Baer radical, transitivity of normality, hyperabelian Lie algebras, abelian subideals.

MSC 2010: 20E15, 20E07, 20E34, 17B05

Not very many relationships between subgroups in groups are transitive. Of the most essential, here you can specify the relations "to be a subgroup", "to be a subnormal subgroup", "to be an ascendant subgroup". At the same time, such important relationships as "to be a normal subgroup", "to be a subnormal subgroup with fixed defect", "to be a permutable subgroup", "to be a pronormal subgroup", "to be an abnormal subgroup" no longer possess the transitivity property. Moreover, the transitivity property turned out to be a rather strict restriction, which in many cases made it possible to obtain an observable description of the corresponding groups ( see [GW1957, AK1958, RD1964, ZG1964, RD1968, RD1973, BBR1999, KS2002, KS2005, BBKP2007, BBKOP2010]).

If $G$ is a group, in which the relation "to be a normal subgroup" is transitive, then every subnormal subgroup of $G$ is normal. This simple circumstance leads us to problems of the following type: to study groups in which not all, but only subnormal subgroups with some natural fixed property are normal. For example, the paper [GF1985] has considered the groups, whose infinite subnormal subgroups are normal; the paper [FGMS2013] has considered the groups, whose subnormal subgroups of infinite special rank are normal.

I.N. Abramovskii in [AI1973] initiated the study of groups in which the transitivity condition is imposed on only abelian normal subgroups. Actually, he studied the groups with the transitivity of normality for Dedekind normal subgroups, but since he considered such locally finite groups with abelian Sylow p-subgroups, those Dedekind subgroups are abelian. It turns out that the class of such groups in which this transitivity inherited by subgroups coincides with the class of locally finite groups with the transitivity for all normal subgroups.

It is clear that the study of groups with all abelian subnormal subgroups are normal makes sense in those classes of groups in which there exist nontrivial subnormal abelian subgroups. One of these classes is the class of hyperabelian groups. Recall that a group $G$ is called *hyperabelian* if $G$ has an ascending series of normal subgroups, whose factors are abelian.

Clearly that in the groups in which the relation "to be a normal subgroup" is transitive ( such groups are called *T − groups* ) every subnormal abelian subgroup is normal. The soluble T – groups were studied in detail in [GW1957, AK1958, RD1964, RD1968, RD1973]. In particular, it turned out that in finite solvable groups the property "to be a T - group" is inherited by subgroups. This, however, is no longer the case for infinite groups. A group G is called $\bar{T}$ – **group**, if every subgroup of G is a T – group. The description of soluble $\bar{T}$ – groups one can find in the paper [RD1964]. The following type of groups is connected with $\bar{T}$ – groups.

A subgroup H of a group G is called **transitively normal** in G, if H is normal in every subgroup S such that H is subnormal in S [KS2006]. It is not hard to see that a group G is a $\bar{T}$ – group if and only if every subgroup of G is transitively normal in G. The question about the groups whose abelian subgroups are transitively normal naturally raised. From the results of the paper [KO2013] it is possible to obtain that a soluble group G is a $\bar{T}$ – group if and only if every abelian subgroup of G is transitively normal in G. In contrast, the situation with groups whose subnormal Abelian subgroups are normal differs significantly from the situation with T - groups. These differences can already be seen in the following example. Let $\{ p_n \mid n \in \mathbb{N} \}$ be the set of odd primes and let $C_n = <c_n>$ be the cyclic group of order $p_n^n$, $n \in \mathbb{N}$. A group $C_n$ has an automorphism $x_n$, having order $p_n^{n-1}$, $n \in \mathbb{N}$, so that we can consider a natural semidirect product $C_n \leftthreetimes <x_n>$. In the Cartesian product $\mathbf{Cr}_{n \in \mathbb{N}} (C_n \leftthreetimes <x_n>)$ we choose the subgroup $G = T \leftthreetimes <x>$ where $T = \mathbf{Dr}_{n \in \mathbb{N}} C_n$ and $x = (x_n)_{n \in \mathbb{N}}$. It is not hard to prove that every abelian subnormal subgroup of G is normal, but G is not a T – group.

We will need the following concepts.

Recall that group G is called a **Dedekind group** if every subgroup of G is normal. The study of such groups was started by R. Dedekind in the paper [DR1897], and therefore these groups were called Dedekind. Later in the paper [BR1933], R. Baer obtained a complete description of such groups, which is as follows:

*If G is non – abelian Dedekind, then $G = Q \times E \times B$ where Q is a quaternion group, E is an elementary abelian 2 – subgroup and B is a periodic abelian 2′ – subgroup.*

Let G be a group. Denote by $\mathbb{B}(G)$ the subgroup generated by subnormal cyclic subgroups. The subgroup $\mathbb{B}(G)$ is called **the Baer radical** of a group G. The subgroup $\mathbb{B}(G)$ is locally nilpotent and every its finitely generated subgroup is subnormal in G [BR1955[1]]. The group G is called **the Baer group** if $G = \mathbb{B}(G)$.

Let G be a group and B, C be the normal subgroups of G such that $B \leq C$. The factor C/B is called **G – central,** if $C_G(C/B) = G$. The factor C/B is called **G – eccentric,** if $C_G(C/B) \neq G$.

Let G be a group and A be a normal subgroup of G. We construct **the upper G – central series of A** as:

$$<0> = A_0 \leq A_1 \leq \ldots A_\alpha \leq A_{\alpha+1} \leq \ldots A_\gamma$$

where

$A_1 = \zeta_G(A) = \{ a \in A \mid [a, g] = 1 \text{ for all elements } g \in G \}$,
$A_{\alpha+1}/A_\alpha = \zeta_G(A/A_\alpha)$, for all ordinals $\alpha < \gamma$, $A_\lambda = \cup_{\beta < \lambda} A_\beta$ for all limit ordinals $\lambda < \gamma$, and $\zeta_G(A/A_\gamma) = <0>$.
We note that every subgroup of this series is G – invariant.

The last term $A_\gamma$ of this series is called **the upper G – hypercenter of A** and will denoted by $\zeta_G^\infty(A)$.

If $A = A_\gamma$, then A is called **G – hypercentral**; if $\gamma$ is finite, then A called **G – nilpotent**.

A normal subgroup A of G is said to be **G – hypereccentric,** if it has an ascending series

$$<0> \leq C_0 \leq C_1 \leq \ldots C_\alpha \leq C_{\alpha+1} \leq \ldots C_\gamma = C$$

of G – invariant subgroups of A such that each factor $C_{\alpha+1}/C_\alpha$ is a G – eccentric and G – chief, for every $\alpha < \gamma$.

We say that the normal abelian subgroup A of a group G has **the Z(G) – decomposition,** if

$$A = \zeta_G^\infty(A) \oplus \eta_G^\infty(A)$$

where $\eta_G^\infty(A)$ is the maximal G – hypereccentric G – invariant subgroup of A. This concept was introduced by D. I. Zaitsev [ZDI1979]. It is not hard to see that in this case, $\eta_G^\infty(A)$ includes every G – hypereccentric G – invariant subgroup of A, in particular, it is unique.

**THEOREM A.** *Let G be a hyperabelian group, whose subnormal abelian subgroups are normal, B be the Baer radical of G, L be the locally nilpotent radical of G, and R be the locally nilpotent residual of G. Suppose that the Baer radical of G is periodic. Then G satisfies the following conditions:*

(i) G/B is abelian and residually finite, every subgroup of B is G – invariant, in particular B is a Dedekind group;

(ii) $R \leq B$, R is G – hypereccentric, the factor - group G/R is hypercentral, $2 \notin \Pi(R)$ and $\Pi(R) \cap \Pi(L/R) = \emptyset$;

(iii) $C_G(R) = L$;

(iv) if p is an odd prime, then the Sylow p – subgroups of G are nilpotent;

(v) if p is an odd prime, then the Sylow p – subgroup of B is the Sylow p – subgroup of locally nilpotent radical. In particular, the Sylow p – subgroup of locally nilpotent radical of G is abelian.

(vi) **Tor**(L) = R × **Tor**(Z) where Z is the upper hypercenter of G.

(vii) if the Sylow 2 – subgroup of G is nilpotent, then the Sylow 2 – subgroup of L is abelian and L = B; in particular, if the orders of elements of Sylow 2 – subgroup $B_2$ of B are bounded, then the Sylow 2 – subgroups of G are nilpotent;

(viii) if $p \notin \Pi(R)$ and the orders of elements of $B_p$ are bounded, then the hypercenter of G of finite number includes $B_p$; if the orders of elements of $B_p$ are not bounded, then the hypercenter of G of the number $\omega$ includes $B_p$;

(ix) the hypercentral length of G/R is at most $\omega + 1$;

(x) if the Sylow 2 – subgroup D of G is not nilpotent, then D satisfies the following conditions:

(xa) $B \cap D = B_2$ is abelian and orders of elements of $B_2$ are not bounded;

(xb) $D = C <a>$ where $C = C_D(B_2)$, $a^2 \in B$ and $b^a = b^{-1}$ for every element $b \in B_2$;

(xc) $B_2 \cap C_D(a)$ is an elementary abelian 2 – subgroup;

(xd) $[C, C] \leq B_2 \cap C_D(a)$, in particular, $[C, C]$ is an elementary abelian 2 – subgroup;

       (xe)   *D includes a normal subgroup $S \geq B_2$ such that the factor – group $D/S$ is elementary abelian and $S = C_S(a)B_2$;*
       (xi)   *if the Sylow 2 – subgroup $L_2$ of $L$ is not nilpotent, then $L_2$ satisfies the following conditions:*

       (xia)   *the Sylow 2 – subgroup $B_2$ of $B$ is abelian and orders of elements of $B_2$ are not bounded;*
       (xib)   *$L_2 = B_2 <a>$, $a^2 \in B_2$ and $x^a = x^{-1}$ for every element $x \in B_2$;*
       (xic)   *$C_{B_2}(a)$ is an elementary abelian 2 – subgroup;*
       (xid)   *either $a^2 = 1$ or $a^4 = 1$.*

*Conversely, in every such group each subnormal abelian subgroup is normal.*

**COROLLARY A1.** *Let $G$ be a hyperabelian group, whose subnormal abelian subgroups are normal and let $B$ be the Baer radical of $G$. Suppose that the Baer radical of $G$ is periodic. If $G$ is locally nilpotent, then $G$ satisfies the following conditions:*

       (i)   *if $p$ is an odd prime, then the Sylow $p$ – subgroup of $G$ is abelian;*
       (ii)   *if the Sylow 2 – subgroup $D$ of $G$ is not nilpotent, then $D$ satisfies the following conditions:*

       (iia)   *the Sylow 2 – subgroup $B_2$ of $B$ is abelian and orders of elements of $B_2$ are not bounded;*
       (iib)   *$D = B_2 <a>$, $a^2 \in B_2$, and $x^a = x^{-1}$ for every element $x \in B_2$;*
       (iic)   *$C_{B_2}(a)$ is an elementary abelian 2 – subgroup;*
       (iid)   *either $a^2 = 1$ or $a^4 = 1$.*

**COROLLARY A2.** *Let $G$ be a periodic hyperabelian group, whose subnormal abelian subgroups are normal, let $B$ be the Baer radical of $G$, $L$ be the locally nilpotent radical of $G$. Suppose that the set $\Pi(G) = \{p_1, p_2, \ldots, p_k\}$ is finite, $k > 1$, and let $p_1 > p_2 > \ldots > p_k$. Then $G$ satisfies the following conditions:*

       (i)   *$B$ is a Dedekind group, then every subgroup of $B$ is $G$ – invariant, $C_G(B) = B$;*
       (ii)   *the factor – group $G/B$ is abelian and finite;*
       (iii)   *$B$ includes the Sylow $p_1$ – subgroup of $G$, the Sylow $p_j$ – subgroup of $G/B$ has a special rank at most $j - 1$, $2 \leq j \leq k - 1$;*
       (iv)   *if $p_k \neq 2$, then the Sylow $p_k$ – subgroup of $G/B$ has a special rank at most $k - 1$, if $p_k = 2$, then the Sylow 2 – subgroup of $G/B$ has a special rank at most $k$;*
       (v)   *if $R$ is a locally nilpotent residual of $G$, then $G = R \leftthreetimes S$ where $S$ is a hypercentral subgroup, moreover, Sylow $2'$ – subgroup of $S$ is nilpotent;*
       (vi)   *if $S_1$ is a subgroup of $G$ such that $G = R \leftthreetimes S_1$, then the subgroups $S$ and $S_1$ are conjugate.*

**COROLLARY A3.** *Let $G$ be a periodic hyperabelian group, whose subnormal abelian subgroups are normal. Suppose that $\Pi(G) = \{p, q\}$ where $p < q$. If $G$ is not locally nilpotent, then $G$ satisfies the following conditions:*

       (i)   *the Sylow $q$ – subgroup $Q$ of $G$ is abelian, $G = Q \leftthreetimes P$ where $P$ is a Sylow $p$ – subgroup of $G$;*
       (ii)   *the Sylow $p$ – subgroups of $G$ are conjugate.*
       (iii)   *$p$ divides $q - 1$;*

*(iv)  a factor – group  $P/C_P(Q)$  is cyclic and  $C_P(Q) \times Q$  is the Baer radical of  G;*

*(v)  if  $p \neq 2$  and  P  is non – abelian, then the orders of elements of  P  are bounded  and  P  is nilpotent;*

*(v)  if  $p = 2$  and  P  is not Dedekind, then  $P/C_P(Q)$  has order  2  and  P  satisfies the following conditions:*

*(va)  $C_P(Q) = B$  is abelian and orders of elements of  B  are not bounded;*
*(vb)  $P = B < a >$, $a^2 \in B$  and  $x^a = x^{-1}$  for every element  $x \in B$;*
*(vc)  $C_B(a)$  is an elementary abelian  2 – subgroup;*
*(vd)  either  $a^2 = 1$  or  $a^4 = 1$.*

This topic is not specific only to group theory. Lie algebras, whose subideals are ideals, have also been studied (see, [SI1969, VV1985]). In the current work, the following description of hyperabelian Lie algebras, whose abelian subideals are ideals, is obtained.

**THEOREM B.** *Let  L  be a hyperabelian  Lie  algebra. Suppose that  L  is non – abelian. If every abelian subideal of  L  is an ideal, then  $L = A \oplus Fd$  where  A  is an abelian ideal of  L  and  $[d, a] = a$  for all elements  $a \in A$. In particular, every subideal of  L  is an ideal.*

**COROLLARY B1.** *Let  L  be a hyperabelian  Lie  algebra. Then every abelian subideal of  L  is an ideal if and only if every subideal of  L  is an ideal.*

An essential generalization of Lie algebras are Leibniz algebras. We note that the Leibniz algebras, whose subideal are ideals, have been studied in papers [KSS2017, KSY2018]. However, the situation with abelian subalgebras of Leibniz algebras is much more complicated. In contrast to Lie algebras, their classes such as solvable, nilpotent, and those who close to them may be completely unsaturated with abelian subalgebras. In the current paper, we give an example of a nilpotent Leibniz algebra of nilpotency class 2, which has a unique abelian subalgebra, and the dimension of this subalgebra is 1. Its construction shows that it is possible to construct examples of this type with very diverse properties. This shows that abelian subalgebras of Leibniz algebras have a very weak effect on the structure of Leibniz algebras, which cannot be said about the influence of nilpotent subalgebras.

## 1. Groups, whose subnormal abelian subgroups are normal

**We begin with some elementary properties.**

**1.1. LEMMA.** *Let  G  be a group, whose  subnormal abelian subgroups are normal, and  L  be the locally nilpotent radical of  G. Then  L  includes the  G – invariant subgroup  K  such that  every subgroup of  K  is  G – invariant and every subnormal subgroup  R  of  G  such that  $K \leq R$  and  $K \neq R$  is not nilpotent.*

**PROOF.** Let  S  be a subnormal nilpotent subgroup of  L. If  x  is an arbitrary element of S, then the subgroup  $< x >$  is subnormal in  L. Then  $< x >$  is subnormal in  G, so that  $< x >$  is normal in  G. It follows that every subgroup of  S  ( and  S  also ) is normal in G. In particular, S  is Dedekind. It follows that  S  is nilpotent of nilpotency class at most  2. Since it is true for each nilpotent subnormal subgroup of  L, the subgroup  K, generated by all nilpotent subnormal subgroup of  L, is nilpotent of nilpotency class at most  2. By above proved every subgroup of  K  is  G – invariant.

**1.2. COROLLARY.** *Let $G$ be a group, whose subnormal abelian subgroups are normal. Then the Baer radical of $G$ includes every nilpotent subnormal subgroup and every subgroup of $B(G)$ is $G$ – invariant ( in particular, $B(G)$ is a Dedekind group ).*

**1.3. COROLLARY.** *Let $G$ be a group, whose subnormal abelian subgroups are normal. Then the factor – group $G/C_G(B(G))$ is abelian.*

**PROOF.** By Corollary 1.2, every subgroup of $B(G)$ is $G$ – invariant. By Corollary 1.2, every subgroup of $B$ is $G$ – invariant. Then $G/C_G(B)$ is abelian ( see, for example, [SchR1994, Theorem 1.5.1 ]).

**1.4. COROLLARY.** *Let $G$ be a group, whose subnormal abelian subgroups are normal. If $G$ is hyperabelian, then every subgroup of $B(G)$ is $G$ – invariant ( in particular, $B(G)$ is a Dedekind group ) and $B(G)$ includes $C_G(B(G))$.*

**PROOF.** By Corollary 1.2, every subgroup of $B(G)$ is $G$ – invariant. Put $B = B(G)$ and $C = C_G(B(G))$, and suppose that $B$ does not include $C$. Then $CB/B$ is not trivial. Since $G$ is hyperabelian, $CB/B$ includes non – trivial abelian $G$ – invariant subgroup $A/B$. The choice of $A$ shows that $A$ is nilpotent. But **Corollary 1.2** implies that, in this case, $B$ must include $A$, and we obtain a contradiction. This contradiction proves the inclusion $C_G(B(G)) \leq B(G)$.

**1.5. COROLLARY.** *Let $G$ be a group, whose subnormal abelian subgroups are normal, and let $B$ the Baer radical of $G$. If $G$ is hyperabelian, then the factor – group $G/B$ is abelian.*

PROOF. By **Corollary 1.2,** every subgroup of $B$ is $G$ – invariant. Then $G/C_G(B)$ is abelian by **Corollary 1.3**. By **Corollary 1.4,** $C_G(B) \leq B$, and therefore $G/L$ is abelian.

**1.6. COROLLARY.** *Let $G$ be a group, whose subnormal abelian subgroups are normal. If $G$ is hyperabelian, then $G$ is hypercyclic.*

PROOF. By **Corollary 1.2,** every subgroup of $B(G)$ is $G$ – invariant, so that $B(G)$ has an ascending series of $G$ – invariant subgroups

$$\langle 1 \rangle = C_0 \leq C_1 \leq \ldots C_\alpha \leq C_{\alpha+1} \leq \ldots C_\beta = B(G),$$

whose factors are cyclic. Since the factor – group $G/B(G)$ is abelian by **Corollary 1.5**, this series can be extended to an ascending series of normal subgroups

$$\langle 1 \rangle = C_0 \leq C_1 \leq \ldots C_\alpha \leq C_{\alpha+1} \leq \ldots C_\beta = B(G) \leq C_{\beta+1} \leq \ldots C_\lambda \leq C_{\lambda+1} \leq \ldots C_\gamma = G,$$

whose factors are cyclic.

**1.7. COROLLARY.** *Let $G$ be a group, whose subnormal abelian subgroups are normal. If $G$ is hyperabelian, then $G$ is locally supersoluble.*

PROOF. We must only note that every hypercyclic group is locally supersoluble [BR1955, Theorem 1].

**1.8. PROPOSITION.** *Let $G$ be a group, whose subnormal abelian subgroups are normal and let $B$ be the Baer radical of $G$. Suppose that $B$ is not periodic. If $G$ is hyperabelian, then $G$ satisfies the following conditions:*

*(i) B is abelian;*

*(ii)* $G = B \langle a \rangle$ and $x^a = x^{-1}$ for every element $x \in B$;
*(iii)* $C_B(a)$ is an elementary abelian $2$–subgroup;
*(iv)* either $a^2 = 1$ or $a^4 = 1$.

*Conversely, in every such group each subnormal abelian subgroup is normal.*

PROOF. By **Corollary 1.2**, B is a Dedekind group. Being non–periodic, B is abelian. Using **Corollary 1.2** we obtain that every subgroup of B is G–invariant. Then $G/C_G(B)$ has order 2 (see, for example, [SchR1994, Theorem 1.5.7]). **Corollary 1.4** shows also that $C_G(B) = B$, so that $G = B \langle a \rangle$ for some element $a \in G$ such that $a^2 \in B$. Moreover, $x^a = x^{-1}$ for every element $x \in L$ (see, for example, [SchR1994, Theorem 1.5.7]). If $y \in C_B(a)$, then we obtain that $y = y^a = y^{-1}$. It follows that $y^2 = 1$. Suppose that $a^2 \neq 1$. Since $a^2 \in C_L(a)$, by above proved $(a^2)^2 = 1$.

Conversely, suppose that a group G satisfies all above condition. Let A be an arbitrary subnormal abelian subgroup of G. If B includes A, then by (ii) every cyclic subgroup of A is G–invariant. It follows that A is a normal subgroup of G. Suppose now that B does not include A. Since B is a maximal subgroup of G, $G = BA$. But in this case, G is nilpotent [HHKL1995, Lemma 4]. However, the fact that B is not periodic together with condition (ii) show that it is impossible. This contradiction shows that B includes A.

**1.9. LEMMA.** *Let G be a group, whose subnormal abelian subgroups are normal, and let B be the Baer radical of G. Suppose that B is periodic. Then every abelian subgroup of $\zeta(B)$ has the $Z(G)$–decomposition.*

PROOF. Let A be an arbitrary subgroup of $\zeta(B)$. By **Corollary 1.2** A and each its subgroup are G–invariant. By **Corollary 1.3** the factor–group $G/C_G(B)$ is abelian. Let S be an arbitrary finite subgroup of A. The inclusion $C_G(B) \leq C_G(S)$ implies that $G/C_G(S)$ is abelian. Then a subgroup S has the $Z(G)$–decomposition: $S = \zeta_G^\infty(S) \oplus \eta_G^\infty(S)$ [DKS2017, Corollary 1.6.5]. Let K be a finite subgroup of $\zeta(B)$ including S. Again K is G–invariant. As above proved, K also has the $Z(G)$–decomposition $K = \zeta_G^\infty(K) \oplus \eta_G^\infty(K)$. Clearly, $\zeta_G^\infty(S) \leq \zeta_G^\infty(S)$. As above noted, $\eta_G^\infty(K)$ includes every G–hypereccentric G–invariant subgroup of K, in particular, $\eta_G^\infty(S) \leq \eta_G^\infty(K)$. Taking into account these arguments, it is not difficult to show that $A = \zeta_G^\infty(A) \oplus \eta_G^\infty(A)$.

**1.10. COROLLARY.** *Let G be a group, whose subnormal abelian subgroups are normal, and let B be the Baer radical of G. Suppose that B is periodic. If $p \in \Pi(B)$, $p \neq 2$, and P is the Sylow p–subgroup of B, then either the upper hypercenter of G includes P, or P is G–hypereccentric.*

PROOF. Since $p \neq 2$, the above noted structure of a Dedekind group indicates that $P \leq \zeta(B)$. **Lemma 1.9** implies that P has the $Z(G)$–decomposition $P = \zeta_G^\infty(P) \oplus \eta_G^\infty(P)$. Suppose the contrary, let both subgroups $\zeta_G^\infty(P), \eta_G^\infty(P)$ are not trivial. Then we can choose the elements $1 \neq c \in \zeta_G(P), 1 \neq d \in \eta_G^\infty(P)$ such that $|c| = |d| = p$. Since subgroup $\langle d \rangle$ is G–invariant, $G/C_G(\langle d \rangle)$ is a non–trivial cyclic group and its order divides $p - 1$. Let g be an element such that $G/C_G(\langle d \rangle) = \langle gC_G(\langle d \rangle) \rangle$. Then $d^g = d^k$ where $1 < k < p$. We have $c^g = c$. Clearly, $cd \notin \zeta_G(P)$. The subgroup $\langle cd \rangle$ is G–invariant, so that $cd \neq (cd)^g = (cd)^m = c^m d^m$ where $1 < m < p$. On the other hand, $c^m d^m = (cd)^g = c^g d^g = cd^k$. It follows that $m \equiv 1 \pmod{p}$ and $m \equiv k \pmod{p}$, and we obtain a contradiction, which proves the result.

**1.11. LEMMA.** *Let $G$ be a group, whose subnormal abelian subgroups are normal and let $B$ be the Baer radical of $G$. Suppose that $B$ is periodic. Then the upper hypercenter of $G$ includes the Sylow $2$ – subgroup of $B$.*

PROOF. Denote by $D$ the Sylow $2$ – subgroup of $B$. By **Corollary 1.2,** $D$ is a Dedekind group. Suppose first that $D$ is abelian. Let

$$D_k = \Omega_k(D) = \{\, d \mid d \in D \text{ and } |d| \text{ divides } p^k \,\}.$$

Consider a factor $D_{k+1}/D_k$. Since every subgroup of $D$ is $G$ – invariant, by **Corollary 1.2**, every subgroup of $D_{k+1}/D_k$ is also $G$ – invariant. We note that the center of each group includes normal cyclic subgroup of order $2$. Since $D_{k+1}/D_k$ is elementary abelian $2$ – group, we obtain that the factor $D_{k+1}/D_k$ is $G$ – central. Since it is true for every positive integer $k$, the equality $D = \bigcup_{k \in \mathbb{N}} D_k$ shows that the upper hypercenter of $G$ includes $D$.

Suppose that $D$ is not abelian. Then $D = Q \times A$ where $Q$ is a quaternion group and $A$ is an elementary abelian $2$ – subgroup. Then $\zeta(D) = \zeta(Q) \times A$, in particular, $\zeta(D)$ and the factor $D/\zeta(D)$ are elementary abelian $2$ – groups. As it was done above, we can prove that $\zeta(D)$ and $D/\zeta(D)$ are $G$ – central, and we obtain that the upper hypercenter of $G$ includes $D$.

**1.12. PROPOSITION.** *Let $G$ be a group, whose subnormal abelian subgroups are normal, $B$ be the Baer radical of $G$, and $R$ be the locally nilpotent residual of $G$. Suppose that $B$ is periodic. If $G$ is hyperabelian, then $B$ includes $R$, $2 \notin \Pi(R)$, $\Pi(R) \cap \Pi(B/R) = \varnothing$. Furthermore, $B$ is $G$ – hypereccentric and the factor – group $G/B$ is hypercentral.*

PROOF. We have $B = D \times S$ where $D$ is a Sylow $2$ – subgroup of $D$ and $S$ is a Sylow $2'$ – subgroup of $B$. Denote by $\pi$ the set of primes $p \in \Pi(B)$ such that the upper hypercenter of $G$ includes the Sylow $p$ – subgroup of $B$. **Lemma 1.11** shows that $2 \in \pi$. By **Corollary 1.10,** we have $B = A \times S$ where $S$ is a Sylow $\pi$ – subgroup of $B$ and $A$ is a Sylow $\pi'$ – subgroup of $B$. By **Corollary 1.5,** the factor-group $G/B$ is abelian. It follows that the factor – group $G/A$ is hypercentral. It follows that $R \leq A$.

Suppose that $R \neq A$. Then there exists a prime $q \in \Pi(B) \setminus \pi$ such that for the Sylow $q$ – subgroup $Q$ of $B$ we have $R \cap Q \neq Q$. Since every subgroup of $Q$ is $G$ – invariant, every subgroup of $QR/R$ is $G$ – invariant. The fact that $G/R$ is residually locally nilpotent implies that every subgroup of $G/R$, having prime order $q$, is $G$ – central. It follows that $QR/R$ has the non – trivial $G$ – central factors. The $G$ – isomorphism $Q/(R \cap Q) \cong_G QR/R$ shows that $Q/(R \cap Q)$ has non – trivial $G$ – central factors. However the choice of $Q$ and **Corollary 1.10** imply that every $G$ – chief factor of $Q/(R \cap Q)$ must be $G$ – eccentric. This contradiction proves the equality $A = R$.

**1.13. LEMMA.** *Let $G$ be a hyperabelian group, whose subnormal abelian subgroups are normal, and let $B$ be the Baer radical of $G$. Suppose that $B$ is periodic. The following assertions hold:*
   *(i) if $p$ is an odd prime, then the Sylow $p$ – subgroups of $G$ are nilpotent;*
   *(ii) if $p$ is an odd prime, then the Sylow $p$ – subgroup of $B$ is the Sylow $p$ – subgroup of the locally nilpotent radical. In particular, the Sylow $p$ – subgroup of the locally nilpotent radical of $G$ is abelian.*

PROOF. Let $p$ be an odd prime, and $P$ be the Sylow $p$ – subgroup of $G$. **Corollary 1.2** shows that every subgroup of $P \cap B = B_p$ is $G$ – invariant, in particular, $B_p$ is a Dedekind group. The fact that p is odd implies that $B_p$ is abelian. If the orders of elements of $B_p$ are bounded, then the factor – group $P/C_P(B_p)$ is finite cyclic p – group ( see, for example, [SchR1994, Theorem 1.5.6 ]). Let x be an element such that $P = <x> C_P(B_p)$. Since $P/B_p = P/(P \cap B) \cong PB/B$ is abelian, both subgroups $C_P(B_p)$ and $<x> B_p$ are normal in P. The subgroup $C_P(B_p)$ is nilpotent, and the subgroup $<x> B_p$ is also nilpotent ( see, for example [CC2008, Corollary 1.77] ). Being a product of two normal nilpotent subgroup $<x> B_p$ and $C_P(B_p)$, P is nilpotent.

If the orders of elements of $B_p$ are not bounded, then the factor – group $P/C_P(B_p)$ is isomorphic to some subgroup of the multiplicative group of the ring of p – adic integers ( see, for example, [SchR1994, Theorem 1.5.6 ] ). But the last group has no elements of order p ( see, for example, [SJP1970, Chapter II, § 3, Theorem 2 ]). Thus $P = C_P(B_p)$, so that P is nilpotent.

Let now $P_1$ be the Sylow p – subgroup of L. As it was proved above, $P_1$ is nilpotent. Then **Corollary 1.2** implies that B includes $P_1$, in particular, $P_1$ is abelian.

**1.14. LEMMA.** *Let G be a hyperabelian group, whose subnormal abelian subgroups are normal, let B be the Baer radical of G and D be the Sylow 2 – subgroup of G. If D is not nilpotent, then D satisfies the following conditions:*

(i) $B \cap D = B_2$ *is abelian and orders of elements of $B_2$ are not bounded;*

(ii) $D = C <a>$ *where* $C = C_D(B_2)$, $a^2 \in B$, *and* $b^a = b^{-1}$ *for every element* $b \in B_2$;

(iii) $B_2 \cap C_D(a)$ *is an elementary abelian 2 – subgroup;*

(iii) $[C, C] \leq B_2 \cap C_D(a)$; *in particular,* $[C, C]$ *is an elementary abelian 2 – subgroup;*

(iv) *D includes a normal subgroup* $S \geq B_2$ *such that the factor – group D/S is elementary abelian and* $S = C_S(a)B_2$.

PROOF. By **Corollary 1.2,** $B_2$ is a Dedekind group, moreover, every subgroup of $B_2$ is G – invariant. Suppose first that $B_2$ is abelian and the orders of its elements are bounded. Then $D/C_D(B_2)$ is a finite cyclic group (see, for example, [SchR1994, Theorem 1.5.6 ]). Put $C = C_D(B_2)$ and let a be an element such that $D = <a> C$. Since $D/B_2 = D/(B \cap D) \cong DB/B$, the fact that G/B is abelian implies that $D/B_2$ is also abelian. Then the both subgroups C and $<a> B_2$ are normal in D. A subgroup C is clearly nilpotent and a subgroup $<a> B_2$ is nilpotent ( see, for example [CC2008, Corollary 1.77] ). It follows that $D = (<a> B_2) C$ is also nilpotent.

If $B_2$ is not abelian, then as we noted above, $B_2 = Q \times E$ where Q is a quaternion group and E is an elementary abelian 2 – subgroup. Since, by **Corollary 1.2,** every cyclic subgroup of E is G – invariant, the center of G includes $E \times \zeta(Q)$. By the same reason, the center of $G/(E \times \zeta(Q))$ includes $B_2/(E \times \zeta(Q))$. It follows that the second hypercenter of D includes $B_2 = B \cap D$. The isomorphism $D/(B \cap D) \cong DB/B$ and the fact that G/B is abelian imply that $D/B_2$ is abelian. It follows that D is nilpotent. The contradictions, which have been obtained in both above cases, show that $B_2$ must be abelian and the orders of its elements are not bounded.

Since $D/B_2$ is abelian, then we obtain that $D \neq C_D(B_2) = C$. By **Corollary 1.2**, every subgroup of $B_2$ is G – invariant, therefore $D/C = <aC>$ is a cyclic group of order 2 and $b^a = b^{-1}$ for each element $b \in B_2$ ( see, for example, [SchR1994, Theorem 1.5.6 ], [SJP1970, Chapter II, § 3, Theorem 2 ]). If $y \in C_B(a)$, then we obtain that $y = y^a = y^{-1}$. It follows that $y^2 = 1$.

Let $c_1, c_2$ be arbitrary elements of C. The fact that $D/B_2$ is abelian implies that $c_1^a = c_1 b_1, c_2^a = c_2 b_2$ for some elements $b_1, b_2 \in B_2$. We have

$$[c_1, c_2]^a = [c_1{}^a, c_2{}^a] = [c_1b_1, c_2b_2] = b_1{}^{-1}[c_1, c_2b_2]b_1 \ [b_1, c_2b_2] = [c_1, c_2b_2] \ [b_1, c_2b_2] =$$
$$[c_1, b_2] \ b_2{}^{-1}[c_1, c_2]b_2 \ [b_1, b_2] \ b_2{}^{-1}[b_1, c_2]b_2 = [c_1, c_2].$$

Since $D/B_2$ is abelian, $[C, C] \leq B_2$, and using (iii) we obtain that $[C, C]$ is an elementary abelian subgroup.

If $b \in B_2$, then $[a, b] = a^{-1} b^{-1} a b = b^2$, so that $[a, B_2] = B_2{}^2$. Put $E = B_2{}^2$. Put $Z/E = C_{C/E}(aE)$. Since $D/E$ is nilpotent of nilpotency class at most 2, the subgroup $Z/E$ is normal in $D/E$. The fact that $B_2/E$ is an elementary abelian 2 – subgroup implies that $Z$ includes $B_2$. The choice of $Z/E$ shows that $[a, Z] \leq E$. Since $Z$ includes $B_2$ and $[a, B_2] = E$, we obtain the equality $[a, Z] = [a, B_2]$. If $z_1, z_2$ are arbitrary elements of $Z$, then $[a, z_1z_2] = [a, z_2] \ z_2{}^{-1}[a, z_1] \ z_2$. The inclusions $Z \leq C$ and $[a, B_2] \leq B_2$ imply that $z_2{}^{-1}[a, z_1]z_2 = [a, z_1]$, so that $[a, z_1z_2] = [a, z_1] \ [a, z_2]$. Thus we obtain athe equality $[a, Z] = \{ [a, z] \mid z \in Z \}$.

Let now $z$ be an arbitrary element of $Z$. The equality $[a, Z] = [a, B_2]$ implies that there exists an element $b \in B_2$ such that $[a, z] = [a, b]$. It follows that $z^{-1}az = b^{-1}ab$ or $zb^{-1} \in C_Z(a)$. Thus we obtain the equality $Z = C_Z(a) B_2$.

Consider the mapping $\xi_a: C/E \longrightarrow B_2/E$ defined by the rule $\xi_a(xE) = [aE, xE]$, $x \in C$. It is not hard to see that $\xi_a$ is an endomorphism of $C/E$ such that $\text{Ker}(\xi_a) = Z/E$ and $\text{Im}(\xi_a) = [aE, C/E]$. The isomorphism

$$(C/E)/(Z/E) = (C/E)/\text{Ker}(\xi_a) \cong \text{Im}(\xi_a) = [aE, C/E]$$

implies that $C/Z \cong (C/E)/(Z/E)$ is an elementary abelian 2 – group. The equality $D = \langle a \rangle C$ implies that $D/Z = \langle aZ \rangle (C/Z)$. Put $S = \langle a, Z \rangle$. Then $S = \langle a \rangle C_Z(a)B_2 = C_S(a)B_2$. Clearly, the factor – group $D/S$ is elementary abelian.

**1.15. COROLLARY.** *Let $G$ be a hyperabelian group, whose subnormal abelian subgroups are normal, let $B$ be the Baer radical of $G$ and $L_2$ be the Sylow 2 – subgroup of the locally nilpotent radical $L$ of $G$. If $L_2$ is not nilpotent, then $L_2$ satisfies the following conditions:*

*(i) the Sylow 2 – subgroup $B_2$ of $B$ is abelian and orders of elements of $B_2$ are not bounded;*

*(ii) $L_2 = B_2 \langle a \rangle$, $a^2 \in B_2$ and $x^a = x^{-1}$ for every element $x \in B_2$;*

*(iii) $C_{B_2}(a)$ is an elementary abelian 2 – subgroup;*

*(iv) either $a^2 = 1$ or $a^4 = 1$.*

PROOF. If we suppose that the orders of elements of $B_2$ are bounded, then **Lemma 1.13** implies that the Sylow 2 – subgroup $D$ of the group $G$ is nilpotent. In this case, $L_2$ is nilpotent too, and we obtain a contradiction. This contradiction shows that the orders of elements of $B_2$ must be not bounded.

Let $C = C_D(B_2)$. A subgroup $C \cap L_2$ is nilpotent, and using again **Corollary 1.2** we obtain that $B$ includes $C \cap L_2$. Since $L_2$ is not nilpotent, $C$ can not include $L_2$. This means that $C \cap L_2 = B_2$. Thus we can choose an element $a$ such that $a \in L_2$. **Lemma 1.13** shows that (ii) $L_2 = B_2 \langle a \rangle$, $a^2 \in C \cap L_2 = B_2$ and $x^a = x^{-1}$ for every element $x \in B_2$.

**Lemma 1.13** shows also that $C_{B_2}(a)$ is an elementary abelian 2 – subgroup. We have that $a^2 \in B_2$, and since $a^2 \in C_{B_2}(a)$ either $a^2 = 1$ or $a^4 = 1$.

**1.16. LEMMA.** *Let $Q$ be a quaternion group of order 8. Then the group $\text{Pot}(Q)$ of all power automorphisms of $Q$ is an elementary abelian 2 – group of order 4.*

PROOF. Let $Q = \langle a_1, a_2 \rangle$ where $a_1^4 = a_2^4 = 1$, $a_1^2 = a_2^2 = c$, $a_2^{-1}a_1 a_2 = a_1^3$. Let $f$ be an arbitrary power automorphism of $Q$. Then $f(a_1) = a_1^k$, $f(a_2) = a_2^t$, where $t, k \in \{1, 3\}$. Thus for a power automorphism $f$ there is only one situation is possible: $f_1(a_1) = a_1^3$, $f_1(a_2) = a_2$, $f_2(a_1) = a_1$, $f_2(a_2) = a_2^3$, $f_3(a_1) = a_1^3$, $f_3(a_2) = a_2^3$. It is possible to check that each automorphism $f_1, f_2, f_3$ is a power authomorphism, has order 2, and therefore, $\mathbf{Pot}(Q) = \langle f_1, f_2, f_3 \rangle$ is an elementary abelian group of order 4.

PROOF OF THEOREM A.

To prove condition (i), we note that the fact that $G/B$ is abelian follows from **Corollary 1.5,** and the fact that every subgroup of $B$ is $G$ – invariant follows from **Corollary 1.2.** We have $B = \mathbf{Dr}_{p \in \Pi(B)} B_p$ where is the Sylow $p$ – subgroup of $B$ for every prime $p$. Since $B$ is a Dedekind group, $B_p$ is abelian for odd prime $p$. A subgroup $B_2$ is abelian or $B_2 = Q \times E$ where $Q$ is a quaternion group and $E$ is an elementary abelian 2 – subgroup. If orders of elements of $B_p$ are not bounded, then we obtain that $G/C_G(B_p)$ is isomorphic to a subgroup of a multiplicative group of the ring of p – adic integers ( see, for example, [SchR1994, Theorem 1.5.6 ]). This group is residually finite. If the orders of elements of $B_p$ are bounded and $p \neq 2$, then $G/C_G(B_p)$ is a finite cyclic group ( see, for example, [SchR1994, Theorem 1.5.6 ]). If $B_2$ is abelian and the orders of elements of $B_2$ are bounded, then again $G/C_G(B_2)$ is a finite cyclic group ( see, for example, [SchR1994, Theorem 1.5.6 ]). If $B_2$ is not abelian, then as we noted above, $B_2 = Q \times E$ where $Q$ is a quaternion group and $E$ is an elementary abelian 2 – subgroup. Since every cyclic subgroup of $E$ is $G$ – invariant, the center of $G$ includes $E$. It follows that $C_G(B_2) = C_G(Q)$, therefore $G/C_G(B_2)$ is finite. We have $C_G(B) = \cap_{p \in \Pi(B)} C_G(B_p)$. The Remak's theorem shows that the factor – group $G/C_G(B)$ is isomorphic to some subgroup of $\mathbf{Cr}_{p \in \Pi(B)} G/C_G(B_p)$, in particular, $G/C_G(B)$ is residually finite. **Corollary 1.4** shows that $C_G(B) \leq B$. If $B_2$ is abelian, then $B$ is also abelian, so that, in this case, $C_G(B) = B$ and $G/B$ is residually finite. If $B_2$ is not abelian, then $C_G(B) = B \cap C_G(B) = (E \times \zeta(Q)) \times \mathbf{Dr}_{p \in \Pi(B), p \neq 2} B_p$. It follows that $B/C_G(B) = QC_G(B)/C_G(B)$, in particular, $B/C_G(B)$ is finite. We note that finite subgroups of residually finite group are closed in profinite topology. It implies that $G/B$ is residually finite, and (i) is proved,

**Proposition 1.12** implies that $R \leq B$, $2 \notin \Pi(R)$, $R$ is $G$ – hypereccentric and $G/R$ is hypercentral, $\Pi(R) \cap \Pi(B/R) = \varnothing$. The fact that $2 \notin \Pi(R)$ together with **Proposition 1.14** imply that $\Pi(R) \cap \Pi(L/R) = \varnothing$.

Let $g \in C_G(R)$. The choice of $g$ shows that $B \leq \zeta(\langle g, B \rangle)$. Since $\langle g, B \rangle/B$ is hypercentral, $\langle g, B \rangle$ is hypercentral. Since the factor – group $G/B$ is abelian, $\langle g, B \rangle$ is a normal subgroup of $G$. Then $\langle g, B \rangle \leq L$, and (iii) is proved.

Assertions (iv) and (v) follows from **Lemma 1.13**.

Let $\pi = \Pi(R)$, $\sigma = \Pi(L) \setminus \Pi(R)$. Since $G$ is hyperabelian, $\Pi(L) = \Pi(B)$. Let $p \in \sigma$, then by **Corollary 1.10** the upper hypercenter $Z$ of $G$ includes $B_p$. If $p$ is an odd prime, then by (v) $B_p$ is the Sylow $p$ – subgroup of $L$. If $B$ does not includes the Sylow 2 – subgroup $L_2$ of $L$, then **Corollary 1.15** implies that $L_2/B_2$ has order 2. Then the centre of $G/B_2$ includes $L_2/B_2$. This means that $Z$ includes $L_2$. In other words, $Z$ includes $\mathbf{Dr}_{p \in \sigma} L_p$, and (vi) is proved.

Let $D$ be the Sylow 2 – subgroup of $L$. Then the Sylow 2 – subgroup $L_2$ of $L$ is nilpotent. **Corollary 1.2** implies that, in this case, $B$ includes $L_2$. By (v), $B$ includes every Sylow $p$ – subgroup of $L$ for all odd primes $p$. Hence $L = B$. If the orders of elements of $B_2$ are bounded, then we can apply the arguments from the first part of the proof of **Lemma 1.13** and obtain that $D$ is nilpotent. Suppose that $B_2$ is not abelian. Then $B_2 = Q \times E$ where $Q$ is a quaternion group and $E$ is an elementary abelian 2 – subgroup. Since every cyclic subgroup of $E$ is $G$ – invariant, the center of $G$ includes $E \times \zeta(Q)$. By the same reason, the center of $G/(E \times \zeta(Q))$ includes $B_2/(E \times \zeta(Q))$. It follows that the second hypercenter of $D$ includes $B_2 = B \cap D$. The fact that $D/B_2$ is abelian implies that $D$ is nilpotent and (vii) is proved.

Let $p$ be an odd prime and $P_1$ be the Sylow $p$ – subgroup of $L$. By (iii) $P_1$ is nilpotent. By **Corollary 1.2** $B$ includes $P_1$, in particular, $P_1$ is abelian, and (iv) is proved.

Suppose now that $p \notin \Pi(R)$. By **Corollary 1.10,** the upper hypercenter of $G$ includes $B_p$. **Corollary 1.2** implies that, in this case, every factor $\Omega_{j+1}(B_p)/\Omega_j(B_p)$ is $G$ – central, thus the series

$$<1> \leq \Omega_1(B_p) \leq \Omega_2(B_p) \leq \ldots \leq \Omega_j(B_p) \leq \Omega_{j+1}(B_p) \leq \ldots \Omega_\omega(B_p) = \cup_{j \in \mathbb{N}} \Omega_j(B_p)$$

is $G$ – central, and (viii) is proved.

By (viii), the hypercenter of $G/R$ of number $\omega$ includes $B/R$. Since $G/B$ is abelian, the hypercentral length of $G/R$ is at most $\omega + 1$.

Assertions (ix) and (x) follows from **Lemma 1.14** and **Corollary 1.15.**

Conversely, suppose that a group $G$ satisfies all above condition. Let $A$ be an arbitrary subnormal abelian subgroup of $G$. If $B$ includes $A$, then by (i) every cyclic subgroup of $A$ is $D$ – invariant. It follows that $A$ is a normal subgroup of $G$.

Suppose now that $B$ does not include $A$. Since $G/B$ is abelian, the product $BA$ is normal subgroup of $G$. We note that the subgroup $BA$ is nilpotent [HHKL1995, Lemma 4]. But, in this case, the Baer radical of $G$ includes $BA$, and we obtain a contradiction. This contradiction proves that $B$ must include $A$.

PROOF OF COROLLARY A1.

Let $p$ be an odd prime and let $P$ be the Sylow $p$ – subgroup of $G$. **Theorem A** implies that $P$ is nilpotent. Using **Corollary 1.2** we obtain that $P$ must be abelian. An assertion (ii) follows from **Theorem A**.

PROOF OF COROLLARY A2.

Assertion (1) follows from **Theotem A**. The fact that $G/B$ is abelian also follows from **Theotem A**. **Corollary 1.7** implies that a group $G$ is locally supersoluble. Then $G$ has a series of normal subgroups

$$<1> = S_0 \leq S_1 \leq S_2 \leq \ldots S_{k-1} \leq S_k = G$$

where $S_1$ is the Sylow $p_1$ – subgroup, $S_2$ is the Sylow $\{p_1, p_2\}$ subgroup, ..., $S_{k-1}$ is the Sylow $\{p_1, p_2, \ldots p_{k-1}\}$ – subgroup of $G$.

Let $B_1 = B \cap S_1$, $C_1 = C_G(B_1)$. By **Corollary 1.2** every subgroup of $B_1$ is $G$ – invariant, in particular, $B_1$ is a Dedekind group. Since $k > 1$, $p_1 \neq 2$. It implies that $B_1$ is abelian. The Sylow $p_1$ – subgroup $S_1 \cap C_1$ of $C_1$ is normal in $G$. The factor $S_1/B_1 = S_1/(B \cap S_1) \cong S_1B/B$ is abelian, so that $S_1 \cap C_1$ is nilpotent. **Corollary 1.2** shows that $B_1$ includes $S_1 \cap C_1$. Then $S_1/B_1 = S_1/(S_1 \cap C_1) \cong S_1C_1/C_1 \leq G/C_1$. If the orders of elements of $B_1$ are at most $p_1^{t_1}$ for some positive integer $t_1$, then the factor – group $G/C_1$ is isomorphic to a subgroup of a multiplicative group of a ring $Z/p_1^{t_1}Z$ ( see, for example, [SchR1994, Theorem 1.5.6 ]). Then $S_1/B_1$ is finite cyclic $p_1$ – group. In this case, the subgroup $S_1$ is nilpotent ( see, for example [CC2008, Corollary 1.77] ). But **Corollary 1.2** shows that, in this case, B must include $S_1$, i.e $S_1 = B_1$, and in particular, $S_1$ is abelian.

If the orders of elements of $B_1$ are not bounded, then the factor – group $G/C_1$ is isomorphic to some subgroup of multiplicative group of a ring of $p_1$ – adic integers ( see, for example, [SchR1994, Theorem 1.5.6 ] ). But the last group has no elements of order $p_1$ ( see, for example, [SJP1970, Chapter II, § 3, Theorem 2 ]). Thus again $S_1 = B_1$. We note that in both cases $G/C_1$ is a cyclic, whose order divides $p_1 - 1$. In particular, $G/C_1$ does not contain the $p_1$ – elements.

Denote by $B_2$ the Sylow $p_2$ – subgroup of B and let $C_2 = C_G(B_2) \cap C_1$. Since G is locally supersoluble, the Sylow $p_2$ – subgroup $Q_2/S_1$ of $C_1/S_1$ is normal. By its choice, $Q_2$ is locally nilpotent, so that $Q_2 = P_2 \times S_1$ where $P_2$ is the Sylow $p_2$ – subgroup of $Q_2$, and hence the Sylow $p_2$ – subgroup of $C_1$. In particular, we obtain that $P_2$ is normal in G. By **Corollary 1.2**, every subgroup of $B_2$ is G – invariant, in particular, $B_2$ is a Dedekind group. Suppose that $p_2 \neq 2$. Then $B_2$ is abelian. If the orders of elements of $B_2$ are at most $p_2^{t_2}$ for some positive integer $t_2$, then the factor $C_1/C_2$ is isomorphic to a subgroup of a multiplicative group of the ring $Z/p_2^{t_2}Z$ ( see, for example, [SchR1994, Theorem 1.5.6 ]). Then $P_2/(P_2 \cap C_2)$ is a finite cyclic $p_2$ – group. Since $P_2/B_2$ is abelian, $P_2 \cap C_2$ is nilpotent. The fact that $P_2 \cap C_2$ is normal in G together with **Corollary 1.2** shows that B includes $P_2 \cap C_2$, i.e $P_2 \cap C_2 = B_2$. It follows that $P_2/B_2$ is a finite cyclic $p_2$ – group. In this case, the subgroup $P_2$ is nilpotent ( see, for example [CC2008, Corollary 1.77] ), and, using again **Corollary 1.2,** we obtain that $P_2 = B_2$.

If the orders of elements of $B_2$ are not bounded, then the factor $C_1/C_2$ is isomorphic to a subgroup of a multiplicative group of the ring of $p_2$ – adic integers ( see, for example, [SchR1994, Theorem 1.5.6 ] ). But the periodic part of the last group is a cyclic group, whose order is $p_2 - 1$ ( see, for example, [SJP1970, Chapter II, § 3, Theorem 2 ]). Thus again $P_2 = B_2$. It follows that $C_1/C_2$ is a cyclic group, whose order divides $p_2 - 1$. In particular, the Sylow $p_2$ – subgroup of $G/C_2$ is cyclic.

Suppose now that $p_2 = 2$. If $B_2$ is abelian and orders of elements of $B_2$ are bounded, then using the above arguments we obtain that the Sylow 2 – subgroup of $C_1$ coincides with $B_2$. It follows that the Sylow 2 – subgroup of $G/C_2$ is cyclic.

Suppose now that $B_2$ is abelian and the orders of its elements are not bounded. Then the factor $C_1/C_2$ is isomorphic to a subgroup of the multiplicative group of the ring of 2 – adic integers ( see, for example, [SchR1994, Theorem 1.5.6 ] ). But the periodic part of the last group has order 2 ( see, for example, [SJP1970, Chapter II, § 3, Theorem 2 ]). Thus there exists an element a such that $C_1 = <a> C_2$ and $b^a = b^{-1}$ for each element $b \in B_2$. It follows that the Sylow 2 – subgroup of $G/C_2$ has special rank at most 2.

Suppose now that $B_2$ is not abelian. Since every subgroup of $B_2$ is G – invariant, $B_2$ is a Dedekind group. Then as we noted above, $B_2 = Q \times E$ where Q is a quaternion group and E is an elementary abelian 2 – subgroup. Since every cyclic subgroup of E is G – invariant, the center of G includes E. It follows that $C_G(B_2) = C_G(Q)$. In this case, **Lemma 1.16** shows that the Sylow 2 – subgroup of $C_1/C_2$ is elementary abelian and has order at most 4. Hence, in this case, the Sylow 2 – subgroup of $C_1$ is an extension of an elementary abelian subgroup by a finite 2 – group. Hence, it is nilpotent ( see, for example [CC2008, Corollary 1.77] ), so that the Sylow 2 – subgroup of $C_1$ coincides with $B_2$. It follows that the Sylow 2 – subgroup of $G/C_2$ is cyclic.

If $\Pi(G) = \{p_1, p_2\}$, then $C_2 = C_G(B) = B$, so that $G/B$ is a finite $p_2$ – group.

If $k > 2$, then we can apply the above arguments and ordinary induction to prove the assertions (ii) – (iv).

Let $R$ be the locally nilpotent residual of $G$. By **Theorem A,** $R$ is $G$ – hypereccentric, $B = R \times Z$ where $\Pi(R) \cap \Pi(Z) = \emptyset$, and $Z$ is a subgroup of the upper hypercenter of $G$. In the factor – group $G/Z$ the subgroup $B/Z$ is $G$ – hypereccentric, and $(G/Z)/(B/Z)$ is finite and abelian. Then $G/Z = B/Z \leftthreetimes S/Z$, and every complement to $B/Z$ in $G/Z$ is conjugate to $S/Z$ [KOS2010, Proposition 2.5]. The equality $B = R \times Z$ shows that $G = R \leftthreetimes S$.

Let $S_1$ be another subgroup such that $G = R \leftthreetimes S_1$. We have $B = R \times (B \cap S_1)$. Let $\pi = \Pi(B)$, then $B \cap S_1$ is a Sylow $\pi'$ – subgroup of $B$. But $B$ is abelian and its Sylow $\pi'$ – subgroup coincides with $Z$. Hence $B \cap S_1 = Z$, in particular, $S_1$ includes $Z$. In the factor – group $G/Z$ we have $G/Z = B/Z \leftthreetimes S_1/Z$. By above noted, $S_1/Z$ and $S/Z$ conjugate. Then the subgroups $S_1$ and $S$ onjugate.

**1.20. COROLLARY.** *Let $G$ be a periodic hyperabelian group, whose subnormal abelian subgroups are normal. Suppose that $\Pi(G) = \{p, q\}$ where $p < q$. If $G$ is not locally nilpotent, then $G$ satisfies the following conditions:*

(i) *the Sylow $q$ – subgroup $Q$ of $G$ is abelian, $G = Q \leftthreetimes P$ where $P$ is a Sylow $p$ – subgroup of $G$;*

(ii) *the Sylow $p$ – subgroups of $G$ conjugate.*

(iii) *$p$ divides $q - 1$;*

(iv) *the factor – group $P/C_P(Q)$ is cyclic, and $C_P(Q) \times Q$ is the Baer radical of $G$;*

(v) *if $p \neq 2$ and $P$ is non – abelian, then the orders of elements of $P$ are bounded, and $P$ is nilpotent;*

(v) *if $p = 2$ and $P$ is not a Dedekind group, then $P/C_P(Q)$ has order 2, and $P$ satisfies the following conditions:*

(va) *$C_P(Q) = B$ is abelian, and orders of elements of $B$ are not bounded;*

(vb) *$P = B \langle a \rangle$, $a^2 \in B$, and $x^a = x^{-1}$ for every element $x \in B$;*

(vc) *$C_B(a)$ is an elementary abelian 2 – subgroup;*

(vd) *either $a^2 = 1$ or $a^4 = 1$.*

PROOF OF COROLLARY A3.

From **Corollary A2** it follows that $Q$ is abelian and the Baer radical of $G$ includes $Q$. **Corollary 1.10** and the fact that $G$ is not locally nilpotent imply that $Q$ must be $G$ – hypereccentric. Since $G/Q$ is a $p$ – group, $Q$ coincides with locally nilpotent residual of $G$. By **Corollary A2**, $G = Q \leftthreetimes P$ and, clearly, in this case, $P$ is the Sylow $p$ – subgroup of $G$. Moreover, **Corollary A2** implies that the Sylow $p$ – subgroups of $G$ conjugate.

If the orders of elements of $Q$ are at most $q^t$ for some positive integer $t$, then the factor – group $G/C_G(Q) \cong P/C_P(Q)$ is isomorphic to a subgroup of the multiplicative group of the ring $\mathbb{Z}/q^t\mathbb{Z}$ ( see, for example, [SchR1994, Theorem 1.5.6 ]). The fact that $q \neq 2$ implies that the last group is cyclic, and $p$ divides $q - 1$. If the orders of elements of $Q$ are not bounded, then the $P/C_P(Q)$ is isomorphic to a subgroup of the multiplicative group of the ring of $q$ – adic integers ( see, for example, [SchR1994, Theorem 1.5.6 ] ). But the periodic part of the last group is a cyclic, having order $q - 1$.

Let B be the Sylow p – subgroup of the Baer radical of G. Since the Baer radical of G is a Dedekind group, then $B \leq C_P(Q)$. If p is odd, then **Theorem A** implies that P is nilpotent. Then $C_P(Q)$ is nilpotent and **Corollary 1.2** implies that the Baer radical of G includes $C_P(Q)$. Thus $B = C_P(Q)$ and $C_P(Q) \times Q$ is the Baer radical of G.

Suppose that p = 2. If P is nilpotent, then $C_P(Q) \times Q$ is the Baer radical of G. If P is not nilpotent, then the result follows from assertion (xi) of **Theorem A.**

## 2. Lie algebras, whose abelian subideals are ideals

It is almost obvious that if L is a Lie algebra, whose subalgebras are ideals, then L is abelian.

**2.1. LEMMA.** *Let L be a Lie algebra, whose abelian subideals are ideals. Then L includes an abelian ideal K such that every subalgebra of K is an ideal of L and every subideal R of L such that $K \leq R$ and $K \neq R$ is not nilpotent.*

PROOF. Let S be a nilpotent subideal of L. If x is an arbitrary element of S, then a subalgebra $<x> = Fx$ is a subideal of S [AS1974, Lemma 3.7]. Then $<x>$ is a subideal of L, so that a subalgebra $<x>$ is an ideal of L. Since it is true for each nilpotent subideal of L, in the subalgebra K generated by all nilpotent subideals of H, every cyclic subalgebra is an ideal of L. It follows that every subalgebra of K ( and K itself ) is an ideal of L. The choice of K shows that K includes every nilpotent subideal of L. By above noted, K is abelian.

Let L be a Lie algebra over a field F, M be non – empty subset of L and H be a subalgebra of L. Put

$$\mathbf{Ann}_H(M) = \{ a \in H \mid [a, M] = <0> \}.$$

The subset $\mathbf{Ann}_H(M)$ is called the ***annihilator*** or the ***centralizer*** of M in subalgebra H.

It is not hard to see that $\mathbf{Ann}_H(M)$ is a subalgebra of L. Moreover, if M is an ideal of L, then $\mathbf{Ann}_L(M)$ is an ideal of L.

**2.2. LEMMA.** *Let L be a Lie algebra over a field F, and let A be an abelian ideal of L. If every subalgebra of A is an ideal of L, then the factor – algebra $L/\mathbf{Ann}_L(A)$ has dimension 1, and for every element $x \in L$ there exists an element $\sigma_x \in F$ such that $[x, a] = \sigma_x a$ for all elements $a \in A$.*

This assertion is ( for example ) a partial case of Lemma 2.2 of the paper [KSY2020].

The Lie algebra L is called ***hyperabelian*** M if L has an ascending series of ideals whose factors are abelian.

### PROOF OF THEOREM B.

PROOF. By **Lemma 2.1,** L has the greatest nilpotent ideal A of L. Being hyperabelian, L includes a non – zero abelian ideal, so that A is non – zero. Let $C = \mathbf{Ann}_L(A)$, and suppose that A does not include C. Then (C + A)/A is non – zero. Since L is hyperabelian, (C + A)/A includes a non – zero abelian ideal S/A. The inclusion S ≤ C shows that S is nilpotent. But **Lemma 2.1** implies that in this case, A must include S, and we obtain a contradiction. This contradiction proves the inclusion $\mathbf{Ann}_L(A) \leq A$. Since A is abelian, $\mathbf{Ann}_L(A) = A$ and every subalgebra of A is an ideal of L. **Lemma 2.2** implies that L/A has dimension 1, so that L = A ⊕ Fb for some element b ∈ L. Using again **Lemma 2.2,** we obtain that there exists an element β ∈ F such that [b, a] = βa for all elements a ∈ A. If we suppose that β = 0, then L is abelian, and we obtain a contradiction. Hence β ≠ 0. Put d = $\beta^{-1}$ b, then

$$[d, a] = [\beta^{-1} b, a] = \beta^{-1} [b, a] = \beta^{-1}\beta a = a \quad \text{for all elements } a \in A.$$

In conclusion, we now give an example of a nilpotent Leibniz algebra of nilpotency class 2, which has a unique abelian subalgebra, and the dimension of this subalgebra is equal to 1.

**2.3. EXAMPLE.** Let n be an arbitrary positive integer and let V be a vector space over a field Q of rational numbers, having dimension n. Let us define on the space V a positive defined bilinear form Φ. In other words, we can choose a basis { $v_1, \ldots, v_n$ } of the space V such that $\Phi(v_k, v_k) = 1$ for all k, $1 \leq k \leq n$, and $\Phi(v_k, v_j) = 0$ whenever k ≠ j. Put L = V ⊕ Qc, and define the operation [,] on L for the elements of the basis by the following rule: [v, c] = [c, v] = [c, c] = 0, and $[v_k, v_j] = \Phi(v_k, v_j)c$, $1 \leq k, j \leq n$, and expand it bilinearly to all elements of L. Put Z = Qc. Then the centre of L includes Z and the factor – algebra L/Z is abelian. It follows that L is a nilpotent Leibniz algebra, moreover, its nilpotency class is 2.

Let x, y be the elements of L. Then

$$x = \xi c + \sum_{1 \leq j \leq n} \alpha_j v_j, \, y = \eta c + \sum_{1 \leq k \leq n} \beta_k v_k, \alpha_1, \ldots, \alpha_n, \beta_1, \ldots, \beta_n \in F.$$

We have
$$[x, y] = [\xi c + \sum_{1 \leq j \leq n} \alpha_j v_j, \eta c + \sum_{1 \leq k \leq n} \beta_k v_k] = \sum_{1 \leq j \leq n} \sum_{1 \leq k \leq n} \alpha_j \beta_k [v_j, v_k] = (\sum_{1 \leq j \leq n} \alpha_j \beta_j) c.$$

It follows that $[x, x] = (\alpha_1^2 + \ldots + \alpha_n^2)c$. In particular, if x ∉ Z, then $(\alpha_1, \ldots, \alpha_n) \neq (0, \ldots, 0)$ and, therefore, [x, x] ≠ 0. Now if S is a subalgebra of L such that Z does not include S, then S contains an element x ∉ Z. Then [x, x] ≠ 0, and hence a subalgebra S is not abelian. This means that Z is the only non-trivial abelian subalgebra of L.

Leonid A. Kurdachenko
Department of Algebra
Oles Honchar Dnipro National University
Gagarin prospect 72, Dnepropetrovsk, Ukraine 49010
e-mail: kurdachenko@hotmail.com

Javier Otal
Departmento de Matemáticas – IUMA

Universidad de Zaragoza
Pedro Cerbuna 12, 50009 Zaragoza, España
e-mail: otal@unizar.es

Igor Ya. Subbotin
Department of Mathematics and Natural Sciences
National University
5245 Pacific Concourse Drive, LA, CA 90045, USA
e-mail: isubboti@nu.edu